\documentclass{article}

\usepackage{arxiv}

\usepackage[utf8]{inputenc} 
\usepackage[T1]{fontenc}    
\usepackage{hyperref}       
\usepackage{url}            
\usepackage{booktabs}       
\usepackage{amsfonts}       
\usepackage{nicefrac}       
\usepackage{microtype}      

\usepackage{graphicx}
\usepackage{amsbsy}
\usepackage{amsmath}
\usepackage{natbib}

\providecommand{\bo}{\mathbf}
\providecommand{\bs}{\boldsymbol}

\providecommand{\bo}{\mathbf}

\providecommand{\tr}{\mbox{tr}}
\DeclareMathOperator*{\argmin}{argmin}

\providecommand{\SM}{\hat{\bs \mu}_\mathrm{SM}}
\providecommand{\HL}{\hat{\bs \mu}_\mathrm{HL}}
\providecommand{\TRSM}{\hat{\bs \mu}_\mathrm{TRSM}}
\providecommand{\TRHL}{\hat{\bs \mu}_\mathrm{TRHL}}
\providecommand{\sm}{\mathrm{SM}}
\providecommand{\hl}{\mathrm{HL}}

\usepackage{xcolor}
\newcommand\klaus[1]{{\color{black}{{#1}}}}
\newcommand\Una[1]{{\color{black}{{#1}}}}
\newcommand\una[1]{{\color{black}{{#1}}}}
\newcommand\hannu[1]{{\color{black}{{#1}}}}

\newtheorem{assumption}{Assumption}

\newtheorem{lemma}{Lemma}

\newtheorem{theorem}{Theorem}

\newtheorem{definition}{Definition}

\def\IP{{\mathbb P}}

\title{The Asymptotic Properties of the One-Sample Spatial Rank Methods}

\author{
  Jyrki M{\"o}tt{\"o}nen \\
  Department of Mathematics and Statistics \\
  University of Helsinki \\
  Finland\\
  \texttt{jyrki.mottonen@helsinki.fi} \\
   \And
  Klaus Nordhausen \\
  Department of Mathematics and Statistics \\
  University of Jyv\"askyl\"a \\
  Finland\\
  \texttt{klaus.k.nordhausen@jyu.fi} \\
  \AND
  Hannu Oja \\
  Department of Mathematics and Statistics \\
  University of Turku \\
  Finland\\
  \texttt{hannu.oja@utu.fi} \\
  \And
  Una Radojicic \\
  Institute of Statistics \& Mathematical\\ Methods in Economics \\
  Vienna University of Technology \\
  Austria \\
  \texttt{una.radojicic@tuwien.ac.at} 
}

\begin{document}

\maketitle

\begin{abstract}
For a set of $p$-variate data points $\bo y_1,\ldots,\bo y_n$, 
there are several versions
of multivariate median and related multivariate sign test proposed and studied in the literature. 
In this paper we consider the asymptotic properties of the multivariate extension of the Hodges-Lehmann (HL) estimator,
the spatial HL-estimator, and the related test statistic.
The asymptotic behavior of the spatial HL-estimator and the related test statistic when $n$ tends to infinity are collected, reviewed, and proved, some for the first time though being used already for a longer time.
We also derive the limiting behavior of the HL-estimator when both the sample size $n$ and the dimension $p$ tend to infinity.
\end{abstract}


\section{Introduction}\label{section:introduction}

For a set of $p$-variate data points $\bo y_1,\ldots, \bo y_n$, there are several versions
of multivariate median and related multivariate sign test proposed and studied in the literature.
For some general reviews, see \citet{Small1990} and \citet{Oja2013}, for example. The so-called spatial median which minimizes the sum $\sum_{i=1}^n \|\bo  y_i - \bs \mu\|$ with
the Euclidean norm 
$\|\bo y\|=({\bo y}^{\top}\bo y)^{1/2}$
has a very long history started by \citet{Gini1929}  and \citet{Haldane1948}.  \citet{Brown1983} has developed many of the properties of the spatial median.
 Taking the gradient of the objective function, one sees that  the multivariate sample spatial median ${\SM}$ solves the
equation
$\sum_{i=1}^n \bo u(\bo y_i-\hat{\bs\mu}) = \bo 0$ with the spatial sign 
$$
\bo u(\bo y) = 
\begin{cases}
\bo y/\|\bo y\|, & \text{when $\bo y\neq \bo 0$},\\
\bo 0, & \text{when $\bo y= \bo 0$}.
\end{cases}
$$
The spatial sign test statistic $\sum_{i=1}^n \bo u(\bo y_i)$  for $H_0 : \klaus{\bs \mu_\sm} = \bo 0$.
was considered by \citet{Mottonen1995}, for example. See also \citet[Chapter 6][]{Oja2010}.

\bigskip 

In this paper, we consider a multivariate extension of the popular Hodges-Lehmann (HL) estimator \citep{HodgesLehmann1963},
the spatial HL-estimator $\klaus{\HL}$,  which minimizes
the sum $\sum_{i<j} \| (\bo y_i+ \bo y_j)/2 -  \bs\mu\|$ and is closely related to the multivariate 
spatial signed-rank test statistic $\sum_{i<j} \bo  u ((\bo y_i+ \bo y_j)/2)$ for $H_0: \klaus{\bs \mu_\hl}= \bo 0$.
See \citet{Chaudhuri1992} and \citet{Mottonen1995} for early studies of these statistics.
Other well-known multivariate extensions are the vector of marginal HL-estimators \citep{PuriSen:1971} and the HL-estimator based on Oja signed-ranks \citep{HettmanspergerMottonenOja1997}. \citet{HallinPaindaveine2002}
constructed multivariate signed-rank test statistics that are based
on standardized spatial signs or Randles' interdirections  \citep{Randles1989} and the ranks of Mahalanobis distances from the origin.
See also \citet{Mottonen2005} for multivariate generalized spatial signed-rank methods.


\citet{Mottonen2010} provided the detailed results on the limiting behavior of the
spatial median and its affine equivariant modification, so called transformation-retransformation estimator,  when $n\to\infty$. In this paper we collect and review in the same way the behavior of the spatial Hodges-Lehmann estimator when $n\to\infty$. Many of the results and auxiliary results can  be collected from \citet{Arcones1998,Bai1990,Chaudhuri1992,Mottonen1995}.


In this article we also study the limiting behavior of the HL-estimator when both  sample size  $n\to\infty$
and the dimension $p\to\infty$.
\citet{Zou2014}  presented an asymptotic  expansion of the spatial
median under elliptical distributions with identity scatter matrix and applied this expansion to a sign-based test for the
sphericity. \citet{Li2022} provide some results for  the asymptotic behavior of sample spatial median under elliptical distributions when $p$ diverges to infinity at the same rate as $n$.
See also \citet{Cheng2018}, \citet{Feng2016a} and \citet{Feng2016b}.


The paper is structured as follows. In Section~\ref{section:one-sample.spatmed} we recall the basic properties of the spatial median and the transformation-retransformation spatial median. In Sections~\ref{section:one-sample.HL} and \ref{section:one-sample.TR.HL} we consider the asymptotic properties of the spatial HL-estimator and the transformation-retransformation HL-estimator. In  Section~\ref{sec:rank_test} the properties of the spatial signed-rank test are reviewed.
In Section~\ref{sec:highdim} we study the limiting behavior of the HL-estimator \una{in high-dimensional regime, i.e.} when both sample size $n$ and dimension $p$ tend to infinity.
The proofs of all the theorems and lemmas are presented in the Appendix.


\section{Review of the spatial median and the transformation-retransformation spatial median}
\label{section:one-sample.spatmed}

In this section we review the properties of the  spatial median and the transformation-retransformation spatial median. The corresponding proofs of the lemmas and theorems are presented in \cite{Mottonen2010}.

Let $\bo y$ be a $p$-variate random vector with cdf $F$ and $p>1$.
The spatial median of $F$ minimizes the objective function
\begin{eqnarray*}
d_1(\bs\mu) = \mathbb{E}\lbrace\|\bo y-\bs\mu\|-\|\bo y\|\rbrace.
\end{eqnarray*}
Note that, as $| \|\bo y-\bs\mu\|-\|\bo y\| |\leq\|\bs\mu\|$, the expectation always exists.
Let ${\bo y}_1,\ldots,{\bo y}_n$ be a random sample from a $p$-variate distribution
$F$ and $\bo Y = (\bo y_1,\ldots,\bo y_n)^{\top}$ be an 
$n\times p$ matrix of the $n$ observation vectors. 
The multivariate sample spatial median \klaus{$\SM = \SM(\bo Y)$} minimizes the objective
function
\begin{eqnarray*}
d_{1n}(\bs\mu)=
\frac{1}{n}\sum_{i=1}^n
\lbrace\|\bo y_i-\bs\mu\|-\|\bo y_i\|\rbrace.
\end{eqnarray*}
We consider next the distribution of the multivariate sample spatial median
$\klaus{\SM}$ under the following assumption:


\begin{assumption}
(a)
The density of $\bo y_i$ is continuous and bounded in an open neighborhood of
the origin.
(b)
The spatial median of
$\bo y_i$ is unique $\klaus{\bs\mu_\sm}=\bo 0$.
\label{assumption:1}
\end{assumption}


We define the following functions
\begin{eqnarray*}
\bo u(\bo y) & = &\|\bo y\|^{-1}\bo y,\ \ \
\bo A(\bo y)=\|\bo y\|^{-1}\left\lbrack\bo I_p-\|\bo y\|^{-2}\bo y\bo y^{\top}\right\rbrack,
\ \ \ \mbox{and}\ \ \
\\
\bo B(\bo y) & = & \bo u(\bo y)\bo u(\bo y)^{\top}
= \|\bo y\|^{-2}\bo y\bo y^{\top},
\end{eqnarray*}
for $\bo y\neq \bo 0$. For $\bo y = \bo 0$ we set $\bo u(\bo 0)=\bo 0$, $\bo A(\bo 0)=\bo 0$
and $\bo B(\bo 0)=\bo 0$. Note that $\bo u(\bo y)=\nabla\|\bo y\|$ is the
spatial sign vector of vector $\bo y$ (the gradient vector of $\|\bo y\|$) and $\bo A(\bo y)=\nabla^2\|\bo y\|$ is the Hessian matrix of $\|\bo y\|$. We also define the sample statistics
\begin{eqnarray*}
\hat{\bo A} 
= \frac{1}{n}\sum_{i=1}^n\bo A(\bo y_i)
\ \ \
\mbox{and}
\ \ \
\hat{\bo B}
= \frac{1}{n}\sum_{i=1}^n\bo B(\bo y_i).
\end{eqnarray*}

\begin{theorem}
Let $\bo y_1,\ldots,\bo y_n$ be iid observations from a distribution satisfying
Assumption \ref{assumption:1}. Then
$$
{\SM}\stackrel{a.s.}\longrightarrow\bo 0.
$$
\label{theorem:0}
\end{theorem}
\begin{theorem}
Let $\bo y_1,\ldots,\bo y_n$ be iid observations from a distribution satisfying
Assumption \ref{assumption:1}.
Then
$$
\sqrt{n} \klaus{\SM} \stackrel{d}\longrightarrow
N_p(\bo 0, \bo A^{-1}\bo B\bo A^{-1}),
$$
where
$\bo A = \mathbb{E}\lbrace\bo A(\bo y_i)\rbrace$
and
$\bo B = \mathbb{E}\lbrace \bo B(\bo y_i)\rbrace$.
\label{theorem:1}
\end{theorem}

\begin{theorem}
Let $\bo y_1,\ldots,\bo y_n$ be a  iid observations from a distribution satisfying
Assumption \ref{assumption:1}.
Then
$$
\hat{\bo A} \stackrel{\IP}\longrightarrow
\bo A
\ \ \
\mbox{and}
\ \ \
\hat{\bo B}  \stackrel{\IP}\longrightarrow
\bo B.
$$
\label{theorem:2}
\end{theorem}

Theorems \ref{theorem:1} and \ref{theorem:2} imply that we can approximate
the covariance matrix of $\klaus{\SM}$ by
$\frac1n\hat{\bo A}^{-1}\hat{\bo B}\hat{\bo A}^{-1}$.

To compute the spatial median ${\SM}$  there are several algorithms proposed in the literature. For example, the algorithms of \cite{Vardi2000}, \cite{Hossjer1995}, \cite{FritzFilzmoserCroux2012} and \cite{KentErConstable2015}. See also \cite{Oja2010} and \cite{MNM}.

\hannu{
The spatial median $\klaus{\SM}$  is rotation
and shift equivariant, that is, 
$$
\klaus{\SM}(\bo Y\bo O^{\top} + \bo 1_n{\bo a}^{\top}) = 
\bo O \klaus{\SM}(\bo Y) + \bo a
$$
for all orthogonal $p\times p$ matrices $\bo O$ and
for all $p\times 1$ vectors $\bo a$. It is not affine equivariant, however, as $
\klaus{\SM}(\bo Y\bo D^{\top})$ and 
$\bo D \klaus{\SM}(\bo Y)$ may be different for $p\times p$ diagonal matrices $\bo D$. }
We get an affine equivariant version of the spatial median by
using a transformation-retransformation (TR) method. See \cite{Chakraborty1998}. The transformation-retransformation method takes advantage of the properties of scatter matrices.

\begin{definition}
\label{def_scatter}
Let $\bo Y = (\bo y_1,\ldots,\bo y_n)^{\top}$ be an $n\times p$ data matrix. A $p\times p$ matrix (a sample statistic) $\bo S(\bo Y)$
is a scatter matrix if 
it is symmetric, non-negative definite and  affine equivariant in the  sense
$$
\bo S(\bo Y{\bo B}^{\top}+ \bo 1_n{\bo a}^{\top}) = \bo B \bo S(\bo Y){\bo B}^{\top}
$$
for all data matrices $\bo Y$, all nonsingular $p\times p$ matrices $\bo B$ and all $p\times 1$ vectors $\bo a$. 
\end{definition}

For some functionals  $\bo S(\bo Y)$
it only holds that
$
\bo S(\bo Y{\bo B}^{\top}+ \bo 1_n{\bo a}^{\top}) \propto \bo B \bo S(\bo Y){\bo B}^{\top}
$
and  $\bo S$ is then often called a shape matrix. In the following  we use  the word scatter also for shape matrices as they both can be similarly used for our purposes. See \cite{Frahm2009} for  a detailed discussion on affine equivariant shape matrices.


\noindent
\hannu{ Let $\bo S$ be a positive definite scatter matrix. Then our transformation matrix is a symmetric matrix ${\bo S}^{-1/2}$ satisfying ${\bo S}^{-1/2}\ \bo S \ {\bo S}^{-1/2}=\bo I_p$  and the retransformation matrix  $\bo S^{1/2}$ is its symmetric inverse. Another possibility is to use the transformation matrix to invariant coordinates and its inverse to transform back, see \cite{TylerCritchleyDumbgenOja:2009} and \cite{Ilmonen2012}. }
The transformation-retransformation procedure is then as follows.
\begin{description}
\item[(1)]
Take any scatter matrix $\bo S=\bo S(\bo Y)$
\item[(2)]
Standardize the data matrix:
${\bo Y}^* = \bo Y({\bo S}^{-1/2})^{\top}$
\item[(3)]
Find the spatial median for the standardized data matrix:
$\klaus{\SM}({\bo Y}^*)$
\item[(4)]
Retransform the estimator: 
$\klaus{\TRSM}(\bo Y) = {\bo S}^{1/2}
\klaus{\SM}({\bo Y}^*).$
\end{description}

\bigskip

It can be easily seen that the transformation-retransformation estimator $\klaus{\TRSM}(\bo Y)$ is affine equivariant, i.e.
$$
\klaus{\TRSM}(\bo Y\bo B^{\top} + \bo 1_n{\bo a}^{\top}) = 
\bo B \klaus{\TRSM}(\bo Y) + \bo a
$$
for all nonsingular $p\times p$ matrices $\bo B$ and all $p\times 1$ vectors $\bo a$. 

For the transformation-retransformation spatial
median a scatter matrix is needed. One possible choice is the 
Tyler's shape matrix.
Tyler's shape matrix is explained in detail
in \cite{Tyler1987} and reviewed in \citet{TaskinenFrahmNordhausenOja2023}.
A problem with  Tyler's shape matrix is that 
it requires a location value in order to center the data.
For a joint estimation of the spatial median and Tyler's shape
matrix one can use the Hettmansperger-Randles  (HR) algorithm \citep{Hettmansperger2002}. 
We assume (without loss of generality) that the population value of $\bo S$ is
$\bo I_p$, and that $\bo S = \bo S(\bo Y)$ is a $\sqrt{n}$-consistent estimator of $\bo I_p$.
The Hettmansperger-Randles estimator of location and scatter are then defined in the following way.

\begin{definition}
Let ${\bs\mu}$ be a $p\times 1$ vector and $\bo S$ a symmetric positive definite
$p \times p$ matrix, and define
${\bo e}_i = {\bo S}^{-1/2}(\bo y_i - {\bs\mu})$, $i=1,\ldots,n$. The Hettmansperger–Randles (HR) estimator of
location and scatter are the values of ${\bs\mu}$ and $\bo S$ which simultaneously satisfy
$$
\frac{1}{n}\sum_{i=1}^n \bo u({\bo e}_i) = \bo 0
\ \ \text{and}\ \ 
\frac{1}{n}\sum_{i=1}^n \bo u({\bo e}_i)\bo u({\bo e}_i)^{\top} = \frac{1}{p}\bo I_p.
$$
\end{definition}

The following asymptotic result was proved by \cite{Mottonen2010}.

\begin{theorem}
Let $\bo Y = (\bo y_1,\ldots,\bo y_n)^{\top}$  be a random sample from a symmetric distribution around zero satisfying Assumption \ref{assumption:1}. Assume also that scatter matrix
$\bo S= \bo S(\bo Y)$ satisfies 
$\sqrt{n}(\bo S-\bo I_p) = O_{\IP}(1)$. Then 
$\sqrt{n}\klaus{\TRSM}(\bo Y)$ and $\sqrt{n}\klaus{\SM}(\bo Y)$
have the same limiting distribution.
\label{theorem:SMTR_convAB}
\end{theorem}

For a detailed comparison of the spatial median and transformation-retransformation spatial median under ellipticity see \citet{MagyarTyler2011}.

\section{The multivariate spatial HL-estimator}
\label{section:one-sample.HL}

In a univariate context, the pseudo-median of a random variable $y$ with cdf $F$ is defined as the median of $(y_1+y_2)/2$, where $y_1$ and $y_2$ are independent copies of $y$. \cite{HodgesLehmann1963} and \cite{Sen1963} suggested independently an estimator for the pseudo-median which is known nowadays as the Hodges-Lehmann estimator. We consider in the following the corresponding concept in a multivariate setting. 


Let $\bo y$ be a $p$-variate random vector with cdf $F$ and $p>1$.
The spatial Hodges-Lehmann location center of $F$
minimizes the objective function
\begin{eqnarray*}
d_2(\bs\mu) = \mathbb{E}\left\lbrace\left\|\frac{\bo y_1+\bo y_2}{2}-\bs\mu\right\|
-\left\|\frac{\bo y_1+\bo y_2}{2}\right\|\right\rbrace,
\end{eqnarray*}
where $\bo y_1$ and $\bo y_2$ are independent 
\una{random vectors} from $F$. 
See \cite{Chaudhuri1992} and \cite{Mottonen1995}.
As in the
spatial median case, the expectation always exists since the expression between
the braces is always bounded.

Let $\bo y_1,\ldots,\bo y_n$ be a random sample from a $p$-variate distribution
$F$. The multivariate spatial Hodges-Lehmann estimator of the location center
${\bs\mu_\hl}$ is defined as the spatial median of pairwise means, i.e. the
Walsh averages,
\begin{eqnarray*}
\bo z_{i,j} = \frac{\bo y_i+\bo y_j}{2},\ \ \ 1\leq i<j\leq n.
\end{eqnarray*}
The sample spatial Hodges-Lehmann estimator ${\HL}$ thus minimizes the criterion
function
\begin{eqnarray*}
d_{2n}(\bs\mu) = {\binom{n}{2}}^{-1}\sum_{i<j}\left\lbrace
\left\|\frac{\bo y_i+\bo y_j}{2}-\bs\mu\right\|
-\left\|\frac{\bo y_i+\bo y_j}{2}\right\|
\right\rbrace
=
{\binom{n}{2}}^{-1}\sum_{i<j}\left\lbrace
\left\|\bo z_{i,j}-\bs\mu\right\|
-\left\|\bo z_{i,j}\right\|
\right\rbrace\una{.}
\end{eqnarray*}

We consider next the distribution of the multivariate spatial Hodges-Lehmann
estimator ${\HL}$ and spatial rank test statistic under mild assumptions. For the asymptotic theory we assume that
\begin{assumption}
(a) The density of $\bo z_{i,j}$ is continuous and bounded.
(b) The spatial median of
$\bo z_{i,j}$ is unique $\klaus{\bs\mu_\hl}=\bo 0$.
\label{assumption:2}
\end{assumption}

\begin{theorem}
\label{theorem:estconv}
Let $\bo z_{i,j}$, $1\leq i<j\leq n$, be observations from a distribution satisfying
Assumption \ref{assumption:2}.
Then
$$
{\HL}\stackrel{a.s.}\longrightarrow\bo 0.
$$
\end{theorem}

\begin{theorem}
\label{theorem:asydist}
Let $\bo z_{i,j}$, $1\leq i<j\leq n$, be observations from a distribution satisfying
Assumption \ref{assumption:2}. Then
$$
\sqrt{n} \klaus{\HL} \stackrel{d}\longrightarrow
N_p(\bo 0, 4\bo A^{-1}\bo B\bo A^{-1}),
$$
where
$$
\bo A =
\mathbb{E}\lbrace\|\bo z_{1,2}\|^{-1}(\bo I_p-\|\bo z_{1,2}\|^{-2}
\bo z_{1,2}\bo z_{1,2}^{\top})\rbrace
\ \ \text{and}\ \ 
\bo B = 
\mathbb{E}\lbrace
(\|\bo z_{1,2}\|^{-1}
\bo z_{1,2})
(\|\bo z_{2,3}\|^{-1}
\bo z_{2,3}^{\top})
\rbrace.
$$
\end{theorem}

Note that the matrices $\bo A$ and $\bo B$ are computed using dependent Walsh averages 
$$
\bo z_{1,2} = \frac{\bo y_1+\bo y_2}{2}\ \ \text{and}\ \ \bo z_{2,3} = \frac{\bo y_2+\bo y_3}{2},
$$
where $\bo y_1$, $\bo y_2$ and $\bo y_3$ are independent copies from $F$.
The proof of Theorem \ref{theorem:asydist} implies that the covariance matrix of $\klaus{\HL}$ can be approximated by
\begin{eqnarray*}
\frac{1}{n}4\hat{\bo A}^{-1}\hat{\bo B}\hat{\bo A}^{-1},
\end{eqnarray*}
where
$
\hat{\bo A}
=
{\binom{n}{2}}^{-1}
\sum_{i<j}
\|\bo z_{i,j}\|^{-1}(\bo I_p-\|\bo z_{i,j}\|^{-2}
\bo z_{i,j}\bo z_{i,j}^{\top})
$
and\\
$
\hat{\bo B} = 
{\binom{n}{3}}^{-1}
\sum_{i<j<k}
(\|\bo z_{i,j}\|^{-1}
\bo z_{i,j})
(\|\bo z_{j,k}\|^{-1}
\bo z_{j,k}^{\top}).
$

\section{Transformation-retransformation HL-estimator}
\label{section:one-sample.TR.HL}

As the spatial median, the spatial HL-estimator given in Section
\ref{section:one-sample.HL} is only shift and rotation equivariant.
Also in this case an affine equivariant version can be found by using
trans\-formation-retransformation method \citep{Chakraborty1998}.
The procedure is as follows.

\begin{description}
\item[(1)]
Take any scatter matrix $\bo S=\bo S(\bo Y)$
\item[(2)]
Standardize the data matrix:
${\bo Y}^* = \bo Y({\bo S}^{-1/2})^{\top}$
\item[(3)]
Find the spatial HL-estimator for the standardized data matrix:
${\HL}({\bo Y}^*)$
\item[(4)]
Retransform the estimator: 
${\TRHL}(\bo Y) = {\bo S}^{1/2}
{\HL}({\bo Y}^*).$
\end{description}

For the transformation-retransformation Hodges-Lehmann estimator the
signed-rank shape matrix seems the most natural transformation
matrix. 
The Hettmansperger-Randles type simultaneous estimators of location and scatter can be defined in the same way as in the spatial median case:

\begin{definition}
Let ${\bs\mu}$ be a $p\times 1$ vector and $\bo S$ a symmetric positive definite
$p \times p$ matrix, and define
${\bo e}_i = {\bo S}^{-1/2}(\bo y_i - {\bs\mu})$, $i=1,\ldots,n$. The simultaneous estimators of
location and scatter are the values of ${\bs\mu}$ and $\bo S$ for which 
$$
\frac{1}{n}\sum_{i=1}^n \hat{\bo q}({\bo e}_i) = \bo 0
\ \ \text{and}\ \ 
\frac{1}{n}\sum_{i=1}^n \hat{\bo q}({\bo e}_i)\hat{\bo q}({\bo e}_i)^{\top} 
\ \propto\ \bo I_p,
$$
where
$\hat{\bo q}(\bo e) = \frac{1}{n}\sum_{i=1}^n\bo u((\bo e_i + \bo e)/2)$
is the estimated signed-rank function.
\end{definition}

\begin{figure}
    \centering
 \includegraphics[width=0.6\textwidth]{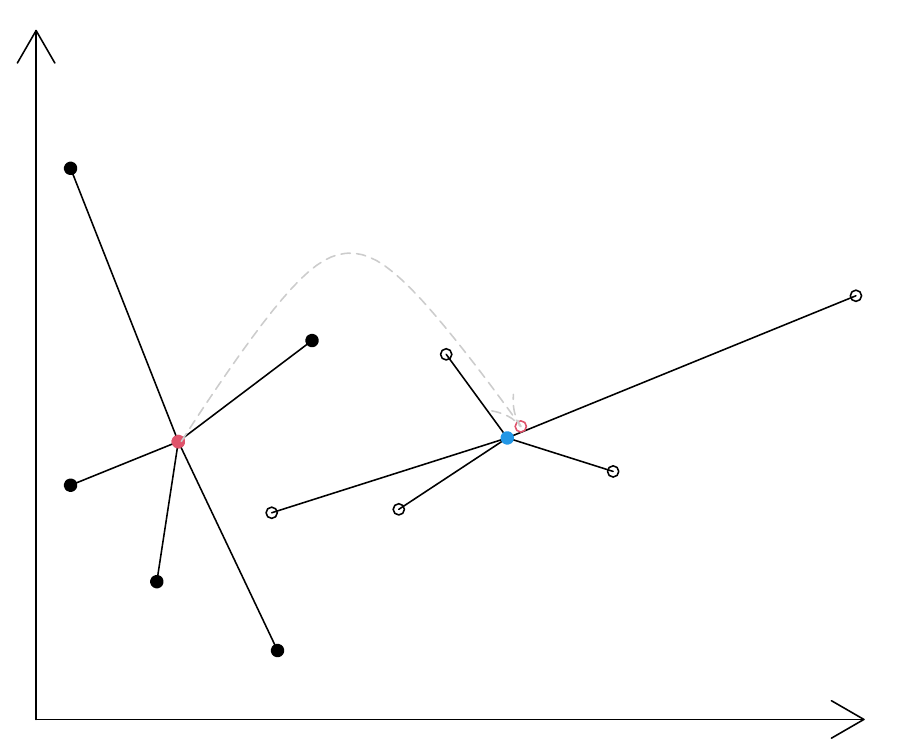}
    \caption{Visualization of effect of the non-affine equivariance of the spatial Hodges-Lehmann estimator. The five solid points have as spatial Hodges-Lehmann estimator the solid red point. Using an affine transformation to these points yields the non-solid points and one can see that the Hodges-Lehmann estimator based on the five transformed points, plotted in blue differs from the transformed original estimate.}
    \label{fig:HL_not_affine_viz}
\end{figure}

\begin{theorem}
Let $\bo Y = (\bo y_1,\ldots,\bo y_n)^{\top}$  be a random sample from a symmetric distribution satisfying Assumption \ref{assumption:2}. Assume also that scatter matrix
$\bo S= \bo S(\bo Y)$ satisfies 
$\sqrt{n}(\bo S-\bo I_p) = O_{\IP}(1)$. Then 
$\sqrt{n} \klaus{\TRHL}(\bo Y)$ and $\sqrt{n} \klaus{\HL} (\bo Y)$
have the same limiting distribution.
\label{theorem:HLTR_convAB}
\end{theorem}

\hannu{
As $\klaus{\HL}$ and  $\klaus{\TRHL}$
are both rotation and shift equivariant, Theorem \ref {theorem:HLTR_convAB}
implies for example that their limiting distributions are the same for
all spherical distributions. For non-spherical elliptical distributions
we expect that the affine equivariant estimator is more efficient even with the small prize needed for the estimation of the scatter matrix. More work is however needed here.
}

To demonstrate the problem of non-affine equivariance of the spatial Hodges-Lehmann estimator  Figure~\ref{fig:HL_not_affine_viz} shows in a very simple case the differences when the five bivariate points are affine transformed. This issue can be avoided by using the transformation-retransformation spatial Hodges-Lehmann estimator.

\section{Spatial signed-rank test}
\label{sec:rank_test}

The spatial signed-rank test statistic for testing the hypothesis \klaus{$H_0: \bs \mu_\hl=\bo 0$} vs. 
\klaus{$H_1: \bs \mu_\hl \neq \bo 0$} can be defined as ($-1$ times) the gradient of the criterion function  (see Section \ref{section:one-sample.HL})
\begin{eqnarray*}
d_{2n}(\bs\mu) = {\binom{n}{2}}^{-1}\sum_{i<j}\left\lbrace
\left\|\frac{\bo y_i+\bo y_j}{2}-\bs\mu\right\|
-\left\|\frac{\bo y_i+\bo y_j}{2}\right\|
\right\rbrace
=
{\binom{n}{2}}^{-1}\sum_{i<j}\left\lbrace
\left\|\bo z_{i,j}-\bs\mu\right\|
-\left\|\bo z_{i,j}\right\|
\right\rbrace
\end{eqnarray*}
evaluated at 
$\bs\mu=\bo 0$:
\begin{eqnarray*}
\bo q_n =
-\frac{\partial 
\una{d_{2n}}(\bs\mu)}{\partial\bs\mu}\bigg|_{\bs\mu=\bo 0} =
{\binom{n}{2}}^{-1}
\sum_{i<j}
\frac{(\bo y_i + \bo y_j)/2}{\|(\bo y_i + \bo y_j)/2\|}
=
{\binom{n}{2}}^{-1}
\sum_{i<j}
\frac{\bo z_{i,j}}{\|\bo z_{i,j}\|}.
\end{eqnarray*}
The general $U$-statistic theory (see e.g. \cite{Hoeffding1948} or \cite{Lee1990}) gives the
following asymptotic result:
\begin{theorem}
Let $\bo z_{i,j}$, $1\leq i<j\leq n$, be observations from a distribution satisfying Assumption~\ref{assumption:2}.
Then
$$
\sqrt{n} \ \bo q_n 
\stackrel{d}\longrightarrow
N_p(\bo 0, 4\bo B),
$$
where
$\bo B = 
\mathbb{E}\lbrace
(\|\bo z_{1,2}\|^{-1}
\bo z_{1,2})
(\|\bo z_{2,3}\|^{-1}
\bo z_{2,3}^{\top})
\rbrace.
$\label{theorem:asytest}
\end{theorem}

Theorem~\ref{theorem:asytest} implies that,  under the null hypothesis $H_0: \klaus{\bs\mu_\hl}=\bo 0$,
$$
((4\bo B)^{-1/2}\sqrt{n}\bo q_n)^{\top}((4\bo B)^{-1/2} \sqrt{n} \bo q_n)
=\frac{n}{4}\bo q_n^{\top}{\bo B}^{-1}\bo q_n \stackrel{d}\longrightarrow \chi_d^2.
$$
If we replace $\bo B$ with the asymptotically consistent estimator $\hat{\bo B}$
we get a spatial signed-rank test statistic
$$
\frac{n}{4}\bo q_n^{\top}{\hat{\bo B}}^{-1}\bo q_n
$$
which is approximately $\chi_d^2$ distributed under null hypothesis $H_0: \klaus{\bs \mu_\hl}=\bo 0$ when $n$ is large.

\begin{figure}
    \centering
 \includegraphics[width=0.49\textwidth]{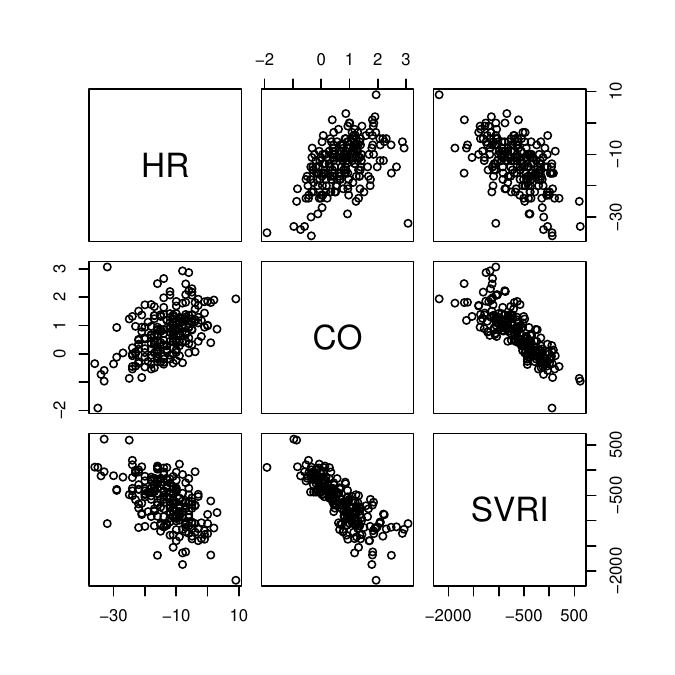}
  \includegraphics[width=0.49\textwidth]{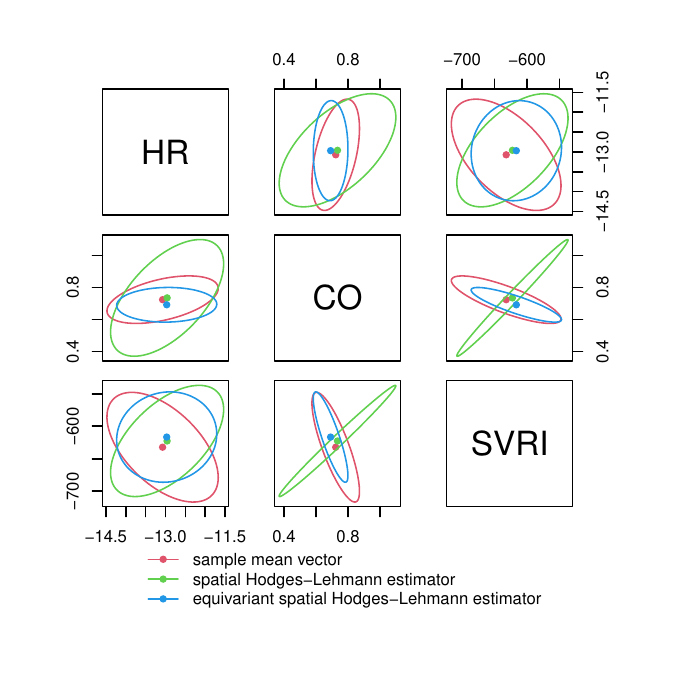}
    \caption{The left panel gives the scatter plot of the differences before and after the tilt for the corresponding variables of the LASERI data. The right panel shows the three estimated locations together with their confidence ellipsoids.}
    \label{fig:laseri}
\end{figure}

To illustrate the spatial signed-rank tests we consider the LASERI data which is publicly available in the R package ICSNP \citep{ICSNP}. For that study several hemodynamic variables were monitored for 233 healthy subjects who were exposed to a passive head-up tilt, i.e. the subjects were first in a lying position, then tilted up before returning again to the lying position. One question here is whether, after the tilt back down and a resting period of 5 minutes, the hemodynamic variables returned to their pre-tilt levels. We will consider the three variables, heart rate (HR), cardiac output (CO), and average systemic vascular resistance index (SVRI), and look at the difference between the average value of the 5th minute before the tilt and the 5th minute after returning to the supine position. The null hypothesis is then that there is no difference, i.e., the location of the difference is $(0,0,0)^\top$. Using the R package MNM \citep{MNM}, we computed for the data the location estimates with confidence ellipsoids as well as the asymptotic test for the transformation-retransformation and non-affine equivariant spatial HL-estimator and as a reference also the mean vector together with Hotelling's $T^2$ test. Figure~\ref{fig:laseri} shows in the left panel the scatter plot of the differences and on the right side the different location estimates with their confidence ellipsoids. The three location estimates are not that different,  however, the non-affine equivariant HL-estimator has \una{rather }
different confidence ellipsoids suffering from the fact that the scales of the three variables are very \una{diverse}
. Comparing all ellipsoids indicates that the affine equivariant HL-estimator is the preferred choice for this data reflecting best the shape of the data and having the smallest volume. However, as for none of the estimates the ellipsoids contain $(0,0,0)^\top$ it comes as no surprise that all three tests yield p-values $<0.0001$.

\section{High-dimensional case}
\label{sec:highdim}
So far the theory above is developed when the dimension $p$ of the distribution is fixed as is usually the case in classical multivariate statistics. In many modern applications however the dimension of the distributions is rather huge leading to a new asymptotic framework known as the high-dimensional \una{regime}. In such a high-dimensional framework the
asymptotic behavior of the sample spatial median was recently 
studied in \citet{Zou2014} assuming spherical symmetry of an underlying distribution, where an asymptotic  expansion of the spatial
median was obtained and further applied to a sign-based test for the sphericity. \cite{Cheng2018} and \cite{Li2022} studied the asymptotic behavior of the sample spatial median under the assumption of elliptical symmetry. Refining the asymptotic representation of the spatial median proposed in \cite{Cheng2018}, \cite{Li2022} gave a modified approximation of the spatial median, where the improvement is in that
the Euclidean norm of the error term reduces to $o_{\IP}(n^{-1/2})$. Such approximation enabled \cite{Li2022} to  establish a central limit theorem for the Euclidean distance of the sample spatial median to the population counterpart, further allowing for the development of one- and
two-sample tests for high-dimensional mean vectors based on sample spatial medians. 
Motivated by the recent developments discussed above, in this section we consider an asymptotic representation of the spatial Hodges-Lehmann estimator in the high dimensional regime and under the assumption of sphericity. In the following, we adopt the notation
$$
\textbf{u}_{i,j}=\textbf{z}_{i,j}/\|\textbf{z}_{i,j}\|,\quad  {r}_{i,j}=\|\textbf{z}_{i,j}\|,\quad c_k=\mathbb{E}(r_{1,2}^{-k}), \,k\geq 1.
$$
Assumption 3(a) we impose further in this section is comparable to those of \cite{Cheng2018}, while Assumption 3(b) is a mild technical assumption ensuring that, in the limit, distribution of $\textbf{z}_{1,2}$ is not fully degenerate, in the sense that the rank of the corresponding covariance is uniformly larger than $1$.
\begin{assumption}
    (a) $\limsup_p\mathbb{E}(r_{1,2}^{-4})<\infty$ and, for $1\leq k\leq 4$, $\limsup_p\mathbb{E}(r_{1,2}^{-k})/ \mathbb{E}(r_{1,2}^{-1})^k=d_k<\infty$.\
    (b) Leading eigenvalue $\lambda_\mathrm{max}(\mathrm{Cov}(\textbf{u}_{1,2}))$ of $\mathrm{Cov}(\textbf{u}_{1,2})$ is uniformly smaller than $1$, i.e. $\limsup_p\lambda_\mathrm{max}(\mathrm{Cov}(\textbf{u}_{1,2}))<1$.
\label{assumption:HD1}
\end{assumption}
Prior to stating the main result, we show that, for general underlying distribution satisfying Assumption~\ref{assumption:HD1}, scaled sample Hodges-Lehman estimator $ \klaus{\HL} $ is bounded in probability.
\begin{theorem}
Let $\textbf{y}_1\dots,\textbf{y}_n$ be a sample from $p$-variate distribution with $\mathbb{E}(\textbf{u}_{1,2})=\bo 0$, and let $n,p\to\infty$. 
Then, under Assumption~\ref{assumption:HD1},
$\sqrt{n} \ a_p \klaus{\HL}  =O_{\IP}(1)$, with choice of $a_p=\mathbb{E}(r_{1,2}^{-1})$. \label{theorem:HD_HL1}
\end{theorem}
Assume now  that $\bo{y}_1,\dots,\bo{y}_n$ is a random sample from an elliptical distribution $\mathcal{E}(\boldsymbol{\mu},2\boldsymbol{\Sigma})$ with mean ${\boldsymbol{\mu}}=\textbf{0}$ and covariance matrix $2\boldsymbol{\Sigma}\propto\textbf{I}_p$. Walsh averages $\textbf{z}_{i,j}=(\bo{y}_{i}+\bo{y}_j)/2$, $i,j=1\dots,n$ then also follow the spherical distribution with location $\boldsymbol{\mu}$ and covariance matrix $\boldsymbol{\Sigma}$, thus admitting representation $\bo z_{i,j}=r_{i,j}\textbf{u}_{i,j}$, where 
$\textbf{u}_{i,j}$ is a random direction distributed uniformly on a unit sphere and independent of $r_{i,j}$. Using this notation it follows that $\boldsymbol{\Sigma}=p^{-1} \mathbb{E}(r_{1,2}^2) \textbf{I}_p$. \Una{Note that Assumption \ref{assumption:HD1} (b) is satisfied for the class of spherical distributions.}
 
\begin{theorem}
Let $\textbf{y}_1\dots,\textbf{y}_n$ be a sample from $p$-variate symmetric distribution around the origin, and let $n,p\to\infty$. Then, under  Assumption~\ref{assumption:HD1}, for the sample Hodges-Lehman estimator $\klaus{\HL}$, \hannu{ $a_p\sqrt{n} \ \klaus{\HL}$ admits the following asymptotic representation:
$$
a_p\sqrt{n} \ \klaus{\HL}
=\sqrt{n}\binom{n}{2}^{-1}\sum_{i<j}\textbf{u}_{i,j}+o_{\IP}(1),
$$
with choice of $a_p=\mathbb{E}(r_{1,2}^{-1})$.}\label{theorem:HD_HL2}
\end{theorem}
Theorem~\ref{theorem:HD_HL2} can be a starting point \una{for developing one- and two-sample tests for the high-dimensional mean which is, however, beyond the scope of this paper.} 
To illustrate Theorem ~\ref{theorem:HD_HL2} however we
performed a small simulation study to demonstrate that $\bs\Delta = \sqrt n\left( \binom{n}{2}^{-1}\sum_{i<j}r_{i,j}^{-1}\right) \klaus{\HL}-\sqrt{n}\binom{n}{2}^{-1}\sum_{i<j}\textbf{u}_{i,j}$ decreases when $n$ and $p$ grow. Note that as $a_p$ is unknown, in the calculation of $\bs\Delta$, we estimate it using its consistent estimator $\binom{n}{2}^{-1}\sum_{i<j}r_{i,j}^{-1}$; see the latter part of proof of Theorem~\ref{theorem:HD_HL2} for more insight. 
\begin{figure}
    \centering
 \includegraphics[width=0.99\textwidth]{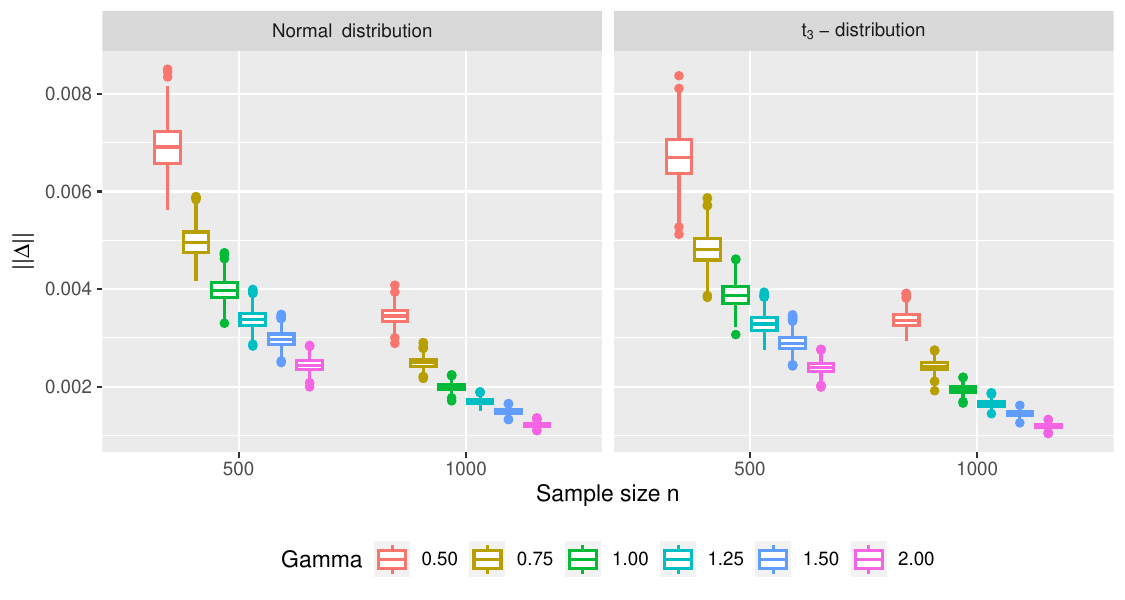}
    \caption{Boxplots of $||\bs\Delta||$ based on 1000 repetitions for different sample sizes $n$ and $\gamma=p/n$ values.}
    \label{fig:HD}
\end{figure}
Based on 1000 repetitions Figure~\ref{fig:HD} shows $||\bs\Delta||$ for different sample sizes $n$ and $\gamma=p/n$'s when $\bo y$ either follows a $p$-variate standard normal distribution or a $p$-variate $t_3$ distribution with $\bs \Sigma = \bo I_p$, and illustrates the limiting behaviour of $\|\bs\Delta\|$, in the high-dimensional regime (Assumption~\ref{assumption:HD1}). More precisely, as $n$ and $p$ increase, the approximation from Theorem~\ref{theorem:HD_HL2} becomes more accurate, i.e. $\|\bs\Delta\|$ decreases, where already for $n=500$ and $p=250$, the median norm of the approximation error is below $0.007$.

\section{Conclusion}
\label{sec:conclusion}
The univariate HL-estimator has a long successful tradition as a location estimator. In this contribution we presented asymptotic results for a multivariate extension based on the concept of spatial signed-ranks yielding the spatial HL-estimator and the transformation-retransformation spatial HL-estimator in the classical multivariate framework where several of the statements appeared earlier in the literature, however often without proofs. What is also novel in this contribution is a first consideration of the spatial HL-estimator in a high-dimensional framework which can be a stepping stone to develop corresponding tests and considerations in the multisample case. \una{Note that in the high-dimensional regime, affine equivariance is not a meaningful concept, and thus we work solely under the framework of orthogonal equivariance.}  
Similarly, recall that under symmetry all the location estimators discussed will estimate the symmetry center while otherwise they are estimating different population quantities.

\section*{Acknowledgments}
The work of KN was partly supported by  HiTEc COST Action (CA21163).
The work of UR is supported by Austrian Science Fund (FWF) (5799-N).

\bibliographystyle{apa-good}
\bibliography{references}

\begin{thebibliography}{41}
\expandafter\ifx\csname natexlab\endcsname\relax\def\natexlab#1{#1}\fi
\expandafter\ifx\csname url\endcsname\relax
  \def\url#1{{\tt #1}}\fi
\expandafter\ifx\csname urlprefix\endcsname\relax\def\urlprefix{URL }\fi

\bibitem[{Arcones(1998)}]{Arcones1998}
Arcones, M.~A. (1998).
\newblock {Asymptotic Theory for $M$-Estimators over a Convex Kernel}.
\newblock {\em Econometric Theory\/}, {\em 14\/}, 387--422.

\bibitem[{Bai et~al.(1990)Bai, Chen, Miao, \& Rao}]{Bai1990}
Bai, Z.~D., Chen, X.~R., Miao, B.~Q., \& Rao, C.~R. (1990).
\newblock {Asymptotic Theory of Least Distances Estimate in Multivariate Linear
  Models}.
\newblock {\em Statistics\/}, {\em 21\/}, 503--519.

\bibitem[{Brown(1983)}]{Brown1983}
Brown, B.~M. (1983).
\newblock {Statistical Uses of the Spatial Median}.
\newblock {\em Journal of the Royal Statistical Society. Series B\/}, {\em
  45\/}, 25--30.

\bibitem[{Chakraborty et~al.(1998)Chakraborty, Chaudhuri, \&
  Oja}]{Chakraborty1998}
Chakraborty, B., Chaudhuri, P., \& Oja, H. (1998).
\newblock {Operating Transformation Retransformation on Spatial Median and
  Angle Test}.
\newblock {\em Statistica Sinica\/}, {\em 8\/}, 767--784.

\bibitem[{Chaudhuri(1992)}]{Chaudhuri1992}
Chaudhuri, P. (1992).
\newblock {Multivariate Location Estimation Using Extension of $R$-Estimates
  Through $U$-Statistics Type Approach}.
\newblock {\em The Annals of Statistics\/}, {\em 20\/}, 897 -- 916.

\bibitem[{Cheng et~al.(2019)Cheng, Liu, Peng, Zhang, \& Zheng}]{Cheng2018}
Cheng, G., Liu, B., Peng, L., Zhang, B., \& Zheng, S. (2019).
\newblock {Testing the Equality of Two High-Dimensional Spatial Sign Covariance
  Matrices}.
\newblock {\em Scandinavian Journal of Statistics\/}, {\em 46\/}, 257--271.

\bibitem[{Davis et~al.(1992)Davis, Knight, \& Liu}]{DAVIS1992}
Davis, R.~A., Knight, K., \& Liu, J. (1992).
\newblock {M-Estimation for Autoregressions with Infinite Variance}.
\newblock {\em {Stochastic Processes and their Applications}\/}, {\em 40\/},
  145--180.

\bibitem[{Feng \& Sun(2016)}]{Feng2016a}
Feng, L., \& Sun, F. (2016).
\newblock {Spatial-Sign Based High-Dimensional Location Test}.
\newblock {\em Electronic Journal of Statistics\/}, {\em 10\/}, 2420 -- 2434.

\bibitem[{Feng et~al.(2016)Feng, Zou, \& Wang}]{Feng2016b}
Feng, L., Zou, C., \& Wang, Z. (2016).
\newblock {Multivariate-Sign-Based High-Dimensional Tests for the Two-Sample
  Location Problem}.
\newblock {\em Journal of the American Statistical Association\/}, {\em 111\/},
  721--735.

\bibitem[{Frahm(2009)}]{Frahm2009}
Frahm, G. (2009).
\newblock {Asymptotic Distributions of Robust Shape Matrices and Scales}.
\newblock {\em Journal of Multivariate Analysis\/}, {\em 100\/}, 1329--1337.

\bibitem[{Fritz et~al.(2012)Fritz, Filzmoser, \&
  Croux}]{FritzFilzmoserCroux2012}
Fritz, H., Filzmoser, P., \& Croux, C. (2012).
\newblock {A Comparison of Algorithms for the Multivariate {L1}-Median}.
\newblock {\em Computational Statistics\/}, {\em 27\/}, 393–410.

\bibitem[{Gini \& Galvani(1929)}]{Gini1929}
Gini, C., \& Galvani, L. (1929).
\newblock Di talune estensioni dei concetti di media ai caratteri qualitativi.
\newblock {\em Metron\/}, {\em 8\/}, 3--209.

\bibitem[{Haldane(1948)}]{Haldane1948}
Haldane, J. B.~S. (1948).
\newblock {Note on the Median of a Multivariate Distribution}.
\newblock {\em Biometrika\/}, {\em 35\/}, 414--417.

\bibitem[{Hallin \& Paindaveine(2002)}]{HallinPaindaveine2002}
Hallin, M., \& Paindaveine, D. (2002).
\newblock {Optimal Tests for Multivariate Location Based on Interdirections and
  Pseudo-Mahalanobis Ranks}.
\newblock {\em The Annals of Statistics\/}, {\em 30\/}, 1103 -- 1133.

\bibitem[{Hettmansperger et~al.(1997)Hettmansperger, M{\"o}tt{\"o}nen, \&
  Oja}]{HettmanspergerMottonenOja1997}
Hettmansperger, T.~P., M{\"o}tt{\"o}nen, J., \& Oja, H. (1997).
\newblock {Affine-Invariant Multivariate One-Sample Signed-Rank Tests}.
\newblock {\em Journal of the American Statistical Association\/}, {\em 92\/},
  1591--1600.

\bibitem[{Hettmansperger \& Randles(2002)}]{Hettmansperger2002}
Hettmansperger, T.~P., \& Randles, R.~H. (2002).
\newblock {A Practical Affine Equivariant Multivariate Median}.
\newblock {\em Biometrika\/}, {\em 89\/}, 851--860.

\bibitem[{Hodges \& Lehmann(1963)}]{HodgesLehmann1963}
Hodges, J.~L., \& Lehmann, E.~L. (1963).
\newblock {Estimates of Location Based on Rank Tests}.
\newblock {\em The Annals of Mathematical Statistics\/}, {\em 34\/}, 598 --
  611.

\bibitem[{Hoeffding(1948)}]{Hoeffding1948}
Hoeffding, W. (1948).
\newblock {A Class of Statistics with Asymptotically Normal Distribution}.
\newblock {\em The Annals of Mathematical Statistics\/}, {\em 19\/}, 293--325.

\bibitem[{H{\"o}ssjer \& Croux(1995)}]{Hossjer1995}
H{\"o}ssjer, O., \& Croux, C. (1995).
\newblock {Generalizing Univariate Signed Rank Statistics for Testing and
  Estimating a Multivariate Location Parameter}.
\newblock {\em Journal of Nonparametric Statistics\/}, {\em 4\/}, 293--308.

\bibitem[{Ilmonen et~al.(2012)Ilmonen, Oja, \& Serfling}]{Ilmonen2012}
Ilmonen, P., Oja, H., \& Serfling, R. (2012).
\newblock {On Invariant Coordinate System ({ICS}) Functionals}.
\newblock {\em International Statistical Review\/}, {\em 80\/}, 93--110.

\bibitem[{Kent et~al.(2015)Kent, Er, \& Constable}]{KentErConstable2015}
Kent, J.~T., Er, F., \& Constable, P. D.~L. (2015).
\newblock {Algorithms for the Spatial Median}.
\newblock In K.~Nordhausen, \& S.~Taskinen (Eds.) {\em Modern Nonparametric,
  Robust and Multivariate Methods: Festschrift in Honour of Hannu Oja\/}, (pp.
  205--224). Cham: Springer.

\bibitem[{Lee(1990)}]{Lee1990}
Lee, A.~J. (1990).
\newblock {\em {$U$-Statistics: Theory and Practice}\/}.
\newblock Routledge.

\bibitem[{Li \& Xu(2022)}]{Li2022}
Li, W., \& Xu, Y. (2022).
\newblock {Asymptotic Properties of High-Dimensional Spatial Median in
  Elliptical Distributions with Application}.
\newblock {\em Journal of Multivariate Analysis\/}, {\em 190\/}, 104975.

\bibitem[{Magyar \& Tyler(2011)}]{MagyarTyler2011}
Magyar, A., \& Tyler, D. (2011).
\newblock {The Asymptotic Efficiency of the Spatial Median for Elliptically
  Symmetric Distributions}.
\newblock {\em Sankhya B\/}, {\em 73\/}, 165--192.

\bibitem[{M{\"o}tt{\"o}nen et~al.(2010)M{\"o}tt{\"o}nen, Nordhausen, \&
  Oja}]{Mottonen2010}
M{\"o}tt{\"o}nen, J., Nordhausen, K., \& Oja, H. (2010).
\newblock {Asymptotic Theory of the Spatial Median}.
\newblock In J.~Antoch, M.~Huskova, \& P.~Sen (Eds.) {\em Nonparametrics and
  Robustness in Modern Statistical Inference and Time Series Analysis\/},
  vol.~7, (pp. 182--193). Institute of Mathematical Statistics.

\bibitem[{M{\"o}tt{\"o}nen \& Oja(1995)}]{Mottonen1995}
M{\"o}tt{\"o}nen, J., \& Oja, H. (1995).
\newblock {Multivariate Spatial Sign and Rank Methods}.
\newblock {\em Journal of Nonparametric Statistics\/}, {\em 5\/}, 201--213.

\bibitem[{M{\"o}tt{\"o}nen et~al.(2005)M{\"o}tt{\"o}nen, Oja, \&
  Serfling}]{Mottonen2005}
M{\"o}tt{\"o}nen, J., Oja, H., \& Serfling, R.~J. (2005).
\newblock {Multivariate Generalized Spatial Signed-Rank Methods}.
\newblock {\em Journal of Statistical Research\/}, {\em 39\/}, 19 -- 42.

\bibitem[{Nordhausen \& Oja(2011)}]{MNM}
Nordhausen, K., \& Oja, H. (2011).
\newblock {Multivariate {$L_1$} Methods: The Package {MNM}}.
\newblock {\em Journal of Statistical Software\/}, {\em 43\/}, 1--28.

\bibitem[{Nordhausen et~al.(2018)Nordhausen, Sirki{\"a}, Oja, \& Tyler}]{ICSNP}
Nordhausen, K., Sirki{\"a}, S., Oja, H., \& Tyler, D.~E. (2018).
\newblock {\em {ICSNP: Tools for Multivariate Nonparametrics}\/}.
\newblock R package version 1.1-1.

\bibitem[{Oja(2010)}]{Oja2010}
Oja, H. (2010).
\newblock {\em {Multivariate Nonparametric Methods with R. An Approach Based on
  Spatial Signs and Ranks}\/}.
\newblock Springer.

\bibitem[{Oja(2013)}]{Oja2013}
Oja, H. (2013).
\newblock {Multivariate Median}.
\newblock In C.~Becker, R.~Fried, \& S.~Kuhnt (Eds.) {\em {Robustness and
  Complex Data Structures: Festschrift in Honour of Ursula Gather}\/}, (pp.
  3--15). Berlin: Springer.

\bibitem[{Puri \& Sen(1971)}]{PuriSen:1971}
Puri, M.~L., \& Sen, P.~K. (1971).
\newblock {\em {Nonparametric Methods in Multivariate Analysis}\/}.
\newblock New York, USA: John Wiley \& Sons.

\bibitem[{Randles(1989)}]{Randles1989}
Randles, R.~H. (1989).
\newblock {A Distribution-Free Multivariate Sign Test Based on
  Interdirections}.
\newblock {\em Journal of the American Statistical Association\/}, {\em 84\/},
  1045--1050.

\bibitem[{Rockafellar(1970)}]{Rockafellar+1970}
Rockafellar, R.~T. (1970).
\newblock {\em {Convex Analysis}\/}.
\newblock Princeton: Princeton University Press.

\bibitem[{Sen(1963)}]{Sen1963}
Sen, P.~K. (1963).
\newblock {On the Estimation of Relative Potency in Dilution (-Direct) Assays
  by Distribution-Free Methods}.
\newblock {\em Biometrics\/}, {\em 19\/}, 532--552.

\bibitem[{Small(1990)}]{Small1990}
Small, C.~G. (1990).
\newblock {A Survey of Multidimensional Medians}.
\newblock {\em International Statistical Review\/}, {\em 58\/}, 263--277.

\bibitem[{Taskinen et~al.(2023)Taskinen, Frahm, Nordhausen, \&
  Oja}]{TaskinenFrahmNordhausenOja2023}
Taskinen, S., Frahm, G., Nordhausen, K., \& Oja, H. (2023).
\newblock {A Review of {Tyler}'s Shape Matrix and Its Extensions}.
\newblock In M.~Yi, \& K.~Nordhausen (Eds.) {\em {Robust and Multivariate
  Statistical Methods: Festschrift in Honor of David E. Tyler}\/}, (pp.
  23--41). Cham: Springer.

\bibitem[{Tyler(1987)}]{Tyler1987}
Tyler, D.~E. (1987).
\newblock {A Distribution-Free $M$-Estimator of Multivariate Scatter}.
\newblock {\em The Annals of Statistics\/}, {\em 15\/}, 234 -- 251.

\bibitem[{Tyler et~al.(2009)Tyler, Critchley, D\"umbgen, \&
  Oja}]{TylerCritchleyDumbgenOja:2009}
Tyler, D.~E., Critchley, F., D\"umbgen, L., \& Oja, H. (2009).
\newblock {Invariant Coordinate Selection}.
\newblock {\em Journal of the Royal Statistical Society. Series B\/}, {\em
  71\/}, 549--92.

\bibitem[{Vardi \& Zhang(2000)}]{Vardi2000}
Vardi, Y., \& Zhang, C.-H. (2000).
\newblock {The Multivariate {$L_1$}-Median and Associated Data Depth}.
\newblock {\em Proceedings of the National Academy of Sciences\/}, {\em 97\/},
  1423--1426.

\bibitem[{Zou et~al.(2014)Zou, Peng, Feng, \& Wang}]{Zou2014}
Zou, C., Peng, L., Feng, L., \& Wang, Z. (2014).
\newblock {Multivariate Sign-Based High-Dimensional Tests for Sphericity}.
\newblock {\em Biometrika\/}, {\em 101\/}, 229--236.

\end{thebibliography}


\appendix
\renewcommand{\thelemma}{A.\arabic{lemma}}
\setcounter{lemma}{0}  

\renewcommand{\thetheorem}{A.\arabic{theorem}}
\setcounter{theorem}{0}  

\newenvironment{myproof}[2] {\paragraph{{\bf {\emph Proof of {\emph {#1}} {\emph {#2}} }}}}{\hfill$\square$}

\section*{Appendix A: Proofs and auxiliary results}
\label{sec:appendix}

\begin{lemma}
Let $\textbf{y}_1,\dots,\textbf{y}_n$ be iid random vectors. 
Let $f(\|(\textbf{y}_i+\textbf{y}_j)/2\|)$ be a function for which 
$\mathrm{Var}(f(\|(\textbf{y}_i+\textbf{y}_j)/2\|))<\infty$. Then, 
$$
{\binom{n}{2}}^{-1}\sum_{i<j} f(\|(\textbf{y}_i+\textbf{y}_j)/2\|)
\stackrel{a.s.}\longrightarrow \mathbb{E}(f(\|(\textbf{y}_1+\textbf{y}_2)/2\|)).
$$
\label{lemma:6}
\end{lemma}

\begin{myproof}{Lemma}{\ref{lemma:6}}
See for example \cite{Lee1990}. 
\end{myproof}

Note that if we define
$$
f(\|(\textbf{y}_i+\textbf{y}_j)/2\|) = 
\left\|\frac{\textbf{y}_i+\textbf{y}_j}{2}-\boldsymbol{\mu}\right\|-\left\|\frac{\textbf{y}_i+\textbf{y}_j}{2}\right\|
$$
Lemma \ref{lemma:6} implies that $d_{2n}(\boldsymbol{\mu})\stackrel{a.s.}\longrightarrow
d_2(\boldsymbol{\mu})$ for all $\boldsymbol{\mu}$.\\

The following key result for convex processes can be found in 
Lemma 2.2 in~\citet{DAVIS1992}
and Theorem 1 in~\citet{Arcones1998}. 
\begin{lemma}
Let $g_n(\bs\theta)$, $\bs\theta\in\mathbb{R}^p$, be a sequence of convex stochastic processes, and
let $g(\bs\theta)$ be a convex (limit) process, meaning that the finite dimensional distributions of $g_n(\bs\theta)$ converge to those of $g(\bs\theta)$. Let further $\hat{\bs\theta},\,\hat{\bs\theta}_1,\dots,\hat{\bs\theta}_n$ be a sequence of random vectors satisfying
$$
g(\hat{\bs\theta})=\inf_{\bs\theta} g(\bs\theta),\quad g_n(\hat{\bs\theta}_n)=\inf_{\bs\theta} g_n(\bs\theta),\quad n=1,2,\dots.
$$
If $\hat{\bs\theta}$ is unique with probability 1 then
$\hat{\bs\theta}_n \stackrel{d}\longrightarrow \hat{\bs\theta}$.
\label{lemma:7}
\end{lemma}

The following lemma is the result of the Lemma 19 of \cite{Arcones1998}. See also \cite{Bai1990} and \cite{Oja2010}.

\begin{lemma}
The accuracies of constant, linear and quadratic approximations of the function $\bs\mu \mapsto \|\bo z - \bs\mu\|$ can be given by
\newline
\begin{tabular}{p{8mm}p{100mm}}
(A1) & $| \|\bo z - \bs\mu\|- \|\bo z\| | \leq \|\bs\mu\|,$ \\ 
(A2) & $| \|\bo z - \bs\mu\|- \|\bo z\| + \bo u^{\top}\bs\mu | \leq 2r^{-1}\|\bs\mu\|^2,$ \\ 
(A3) & $| \|\bo z - \bs\mu\|- \|\bo z\| + \bo u^{\top}\bs\mu - \bs\mu^{\top}(2r)^{-1}[\bo I_p -\bo u\bo u^{\top}]\bs\mu| \leq c r^{-1-\delta}\|\bs\mu\|^{2+\delta}$ for all $0<\delta<1$, \\ 
\end{tabular} 
\newline
where $\bo z = r\bo u$, $r=\|\bo z\|$, $\bo u = \|\bo z\|^{-1}\bo z$ and 
the constant $c$ does not depend on $\bo z$ or $\bs\mu$.
\label{lemma:accuracy}
\end{lemma}

\begin{myproof}{Lemma}{\ref{lemma:accuracy}}
 See e.g. \cite{Oja2010}.
\end{myproof}

\begin{myproof}{Theorem}{\ref{theorem:estconv}}
As  $d_{2n}(\boldsymbol{\mu}) $ and 
$d_{2}(\boldsymbol{\mu}) $
are bounded and convex, 
also 
$\sup_ {\| \boldsymbol {\mu} \|\le C } |d_{2n}(\boldsymbol{\mu})-d_{2}(\boldsymbol{\mu})| \stackrel{a.s.}\longrightarrow 0$, for all $C>0$ (Theorem 10.8 in \cite{Rockafellar+1970}). 
Write  $\HL^* = \argmin_{
\| \boldsymbol {\mu} \|\le C}
d_{2n}(\boldsymbol{\mu})$.
Then 
$d_2(\HL^* )\to 0$ and $\HL^*\to \bo 0$ almost surely as
$d_2(\bo 0)=0$, $|d_{2n}({\HL}^* )-d_{2}({\HL}^* )| \stackrel{a.s.}\longrightarrow 0 $, and
$
d_{2n}({\HL}^* )\le d_{2n}(\bo 0 )
\stackrel{a.s.}\longrightarrow  d_{2}(\bo 0 ) \le
d_{2}({\HL}^*).
$
The result then follows as, using any $C'>C$, the sequence of estimators $\argmin_{
\| \boldsymbol {\mu} \|\le C'}
d_{2n}(\boldsymbol{\mu})$ and  $\HL^*$ are the same almost surely
(use again Theorem 10.8 in \cite{Rockafellar+1970}).
\end{myproof}

\begin{myproof}{Theorem}{\ref{theorem:asydist}}
Define the sample statistic
\begin{eqnarray*}
\hat{\bo A}
=
{\binom{n}{2}}^{-1}
\sum_{i<j}
\|\bo z_{i,j}\|^{-1}(\bo I_p-\|\bo z_{i,j}\|^{-2}
\bo z_{i,j}\bo z_{i,j}^{\top})
\end{eqnarray*}
and the corresponding population value
$$
\bo A =
\mathbb{E}(\hat{\bo A}) =
\mathbb{E}\lbrace\|\bo z_{i,j}\|^{-1}(\bo I_p-\|\bo z_{i,j}\|^{-2}
\bo z_{i,j}\bo z_{i,j}^{\top})\rbrace.
$$
Note that 
$
\|\bo z\|^{-1}(\bo I_p-\|\bo z\|^{-2}
\bo z\bo z^{\top})
$
is the Hessian matrix of $\|\bo z - \bs\mu\|$ evaluated at $\bs\mu=\bo0$.

The approximation (A3) of the Lemma \ref{lemma:accuracy} implies that
\begin{eqnarray*}
\left|
\sum_{i<j}\left\lbrace
\|\bo z_{i,j} - \bs\mu\|- \|\bo z_{i,j}\| + \bo u_{i,j}^{\top}\bs\mu - \bs\mu^{\top}(2r_{i,j})^{-1}[\bo I_p -\bo u_{i,j}\bo u_{i,j}^{\top}]\bs\mu
\right\rbrace
\right|
\leq 
c\|\bs\mu\|^{2+\delta}\sum_{i<j}\frac{1}{r_{i,j}^{1+\delta}} 
\end{eqnarray*}
Replacing $\bs\mu$ by $n^{-1/2}\bs\mu$ and multiplying both sides of the inequality by the constant
$
n{\binom{n}{2}}^{-1}
$
gives the approximation
\begin{eqnarray*}
&&
\bigg|
n
{\binom{n}{2}}^{-1}
\sum_{i<j}
\left\lbrace
\|\bo z_{i,j} - n^{-1/2}\bs\mu\|- \|\bo z_{i,j}\| 
\right\rbrace
+ n^{1/2}
{\binom{n}{2}}^{-1}
\sum_{i<j}
\frac{\bo z_{i,j}^{\top}}{\|\bo z_{i,j}\|}
\bs\mu 
\\
&&
- \bs\mu^{\top}
{\binom{n}{2}}^{-1}
\sum_{i<j}\left[
\frac{1}{2\|\bo z_{i,j}\|}
\left[\bo I_p -
\frac{\bo z_{i,j}\bo z_{i,j}^{\top}}{\|\bo z_{i,j}\|^2}
\right]\right]\bs\mu
\bigg|
\\
&&=
\left|
n
d_{2n}(n^{-1/2}\bs\mu)
+ (\sqrt{n}\bo q_n)^{\top}
\bs\mu 
- 
\frac{1}{2}
\bs\mu^{\top}
\hat{\bo A}\bs\mu
\right|
\leq 
n
{\binom{n}{2}}^{-1}
\frac{c\|\bs\mu\|^{2+\delta}}{n^{(2+\delta)/2}}
\sum_{i<j}\frac{1}{r_{i,j}^{1+\delta}},
\end{eqnarray*}
where  (see \cite{Lee1990})
$\hat{\bo A}\stackrel{a.s.}\longrightarrow\bo A$
and
\begin{eqnarray*}
\sqrt{n}\bo q_n =
\sqrt{n}\left[\left\lbrace{\binom{n}{2}}^{-1}
\sum_{i<j}
\frac{\bo z_{i,j}}{\|\bo z_{i,j}\|}\right\rbrace - \bo 0
\right]
\stackrel{d}\longrightarrow
N_p(\bo 0, 4\bo B),
\end{eqnarray*}
where
$\bo B = 
E\lbrace
(\|\bo z_{1,2}\|^{-1}
\bo z_{1,2})
(\|\bo z_{2,3}\|^{-1}
\bo z_{2,3}^{\top})
\rbrace.
$
We now get under Assumption 2 that
\begin{eqnarray*}
n~d_{2n}(n^{-1/2}\bs\mu)
-
\left(
-(\sqrt{n}~\bo q_n)^{\top}
\bs\mu + 
\frac{1}{2}
\bs\mu^{\top}
\bo A
\bs\mu
\right)
\stackrel{\IP}\longrightarrow 0.
\end{eqnarray*}

We can then apply Lemma \ref{lemma:7} with 
$
g_n(\bs\mu) = n~d_{2n}(n^{-1/2}\bs\mu)
$
and
$
g(\bs\mu) = -\bo h^{\top}\bs\mu + \frac{1}{2}\bs\mu^{\top}\bo A\bs\mu,
$
where $\bo h\sim N_p(\bo 0, 4\bo B)$.
Next, taking the gradient of 
$-(\sqrt{n}~\bo q_n)^{\top}
\bs\mu + 
\frac{1}{2}
\bs\mu^{\top}
\bo A
\bs\mu
$
with respect to $\bs\mu$
and setting it to zero gives ($\bo A$ is nonsingular)
$$
\sqrt{n}\klaus{\HL} 
\stackrel{d}\longrightarrow
N_p(\bo 0, 4\bo A^{-1}\bo B\bo A^{-1}).
$$
 \end{myproof}

\begin{lemma}
Let $\bo Y = (\bo y_1,\ldots,\bo y_n)^{\top}$ be a random sample from a symmetric distribution satisfying Assumption \ref{assumption:1}.
Assume also that scatter matrix $\bo S=\bo S(\bo Y)$ satisfies $\sqrt{n} (\bo S^{-1/2}-\bo I_p)=O_{\IP}(1)$.
Then
$$
\sqrt{n}({\bo q}_{n}(\bo Y^*)-{\bo q}_{n}(\bo Y))\stackrel{\IP}\longrightarrow \bo 0,
$$
where 
$$
{\bo q}_{n}(\bo Y) = 
{\binom{n}{2}}^{-1}\sum_{i<j}\bo u((\bo y_i+\bo y_j)/2)
=
{\binom{n}{2}}^{-1}\sum_{i<j}\bo u(\bo z_{i,j}).
$$
\label{lemma:asyHLtest}
\end{lemma}

\begin{myproof}{Lemma}{\ref{lemma:asyHLtest}}
Let $\bs\Delta_n=\sqrt{n} (\bo S^{-1/2}-\bo I_p)$. Then
$
{\bo S}^{-1/2}=\bo I_p+n^{-1/2} \bs\Delta_n,
$
where $\bs\Delta_n$ is  bounded in probability. 
Using the approximation (B2) in \cite{Mottonen2010}
we obtain
\begin{eqnarray*}
&&\sqrt{n}{\binom{n}{2}}^{-1}\sum_{i<j} \bo u(\bo S^{-1/2}\bo z_{i,j})-
\sqrt{n}{\binom{n}{2}}^{-1} \sum_{i<j} \bo u(\bo z_{i,j}) 
= -{\binom{n}{2}}^{-1}\sum_{i<j}\bo g(\bo z_{i,j})+o_{\IP}(1),
\end{eqnarray*}
where $
\bo g(\bo z_{i,j}) =
[\bo I_p-\bo u(\bo z_{i,j})\bo u(\bo z_{i,j})^{\top}]\bs\Delta_n\bo u(\bo z_{i,j})
$
is symmetrically distributed around the origin as multiplying all the observation by $-1$ simply changes its sign only.  
Let $|\bs\Delta_n|$ denote the Frobenius norm $(\tr(\bs\Delta_n^{\top}\bs\Delta_n))^{1/2}$.
Then ${\binom{n}{2}}^{-1}\sum_{i<j}\bo g(\bo z_{i,j})1_{|\bs\Delta_n|\le M} $ converges  in probability to zero for all $M>0$. As $\IP(|\bs\Delta_n|>M)$ can be made arbitrarily close to zero,
\[
\sqrt{n}{\binom{n}{2}}^{-1}\sum_{i<j}  \bo u(\bo S^{-1/2}\bo z_{i,j})-\sqrt{n}{\binom{n}{2}}^{-1}\sum_{i<j}  \bo u(\bo z_{i,j}) \stackrel{\IP}\longrightarrow 0
\]
and the proof follows.
\end{myproof}

\begin{lemma}
Let $\bo Y = (\bo y_1,\ldots,\bo y_n)^{\top}$ be a random sample from a  distribution satisfying Assumption \ref{assumption:1}.
Assume also that scatter matrix $\bo S=\bo S(\bo Y)$ satisfies $\sqrt{n} (\bo S-\bo I_p)=O_{\IP}(1)$.
 Let $\bo Y^*= \bo Y({\bo S}^{-1/2})^{\top}$.
 Then
\[ 
\bo A(\bo Y^*)-\bo A(\bo Y)\stackrel{\IP}\longrightarrow \bo 0
\ \ \mbox{and}\ \
\bo B(\bo Y^*)-\bo B(\bo Y)\stackrel{\IP}\longrightarrow \bo 0.
\]
\label{TRLemma2}
\end{lemma}

\begin{myproof}{Theorem}{\ref{theorem:HLTR_convAB}}
Let
${\bo S}^{-1/2}=\bo I_p+n^{-1/2} \bs\Delta_n$ and 
${\bo z}_{i,j}^* = (\bo I_p+n^{-1/2} \bs\Delta_n)\bo z_{i,j}$, $i,j=1,\ldots,n$.

The approximation (A3) of the Lemma \ref{lemma:accuracy} gives the result
\begin{eqnarray*}
&&\bigg|
n{\binom{n}{2}}^{-1}
\sum_{i<j}
\left\lbrace
\|\bo z_{i,j}^* - n^{-1/2}\bs\mu\|- \|\bo z_{i,j}^*\|
\right\rbrace 
+
\sqrt{n}
{\binom{n}{2}}^{-1}
\sum_{i<j}
\frac{\bo z_{i,j}^{*\top}}{\|\bo z_{i,j}^*\|}
\bs\mu 
\\
&&- 
\bs\mu^{\top}
\bigg\lbrace
{\binom{n}{2}}^{-1}
\sum_{i<j}
\frac{1}{2\|\bo z_{i,j}^*\|}
\bigg[\bo I_p -
\frac{\bo z_{i,j}^*\bo z_{i,j}^{*\top}}{\|\bo z_{i,j}^*\|^2}
\bigg]
\bigg\rbrace
\bs\mu
\bigg|
\\
&&\leq 
\frac{c\|\bs\mu\|^{2+\delta}}{n^{\delta/2}}
{\binom{n}{2}}^{-1}
\sum_{i<j}
\|(\bo I_p+n^{-1/2} \bs\Delta_n)\bo z_{i,j}\|^{-1-\delta}
\stackrel{\IP}\longrightarrow 0.
\end{eqnarray*}
Lemmas \ref{lemma:7}, \ref{lemma:asyHLtest} and
\ref{TRLemma2}  
imply that
$
\sqrt{n} \ \klaus{\HL}(\bo Y^*) 
$
and 
$
\sqrt{n} \ \klaus{\HL}(\bo Y) 
$
have the same limiting distribution. The result then follows from Slutsky's theorem.
\end{myproof}

\begin{myproof}{Theorem}{\ref{theorem:HD_HL1}}
The HL median is a spatial median of the sample $\textbf{z}_{i,j}$, $i,j=1,\dots,n$ and thus minimizes the objective function $\una{d_{2n}}(\boldsymbol{\mu})=\binom{n}{2}^{-1}\sum_{i<j}\|\textbf{z}_{i,j}-\boldsymbol{\mu}\|-\sum_{i<j}\|\textbf{z}_{i,j}\|$. Write now $a_p\sqrt{n}\boldsymbol{\mu}=C\textbf{v}$, where $C>0$ and $\|\textbf{v}\|=1$. The aim is now to show that for $\varepsilon>0$ there exists $C>0$ such that
$$
\liminf_n\mathbb{P}\left(\inf_{\textbf{v}}\una{d_{2n}}(Ca_p^{-1}n^{-1/2}\textbf{v})>0\right)>1-\varepsilon.
$$
Convexity of $\una{d_{2n}}$ then implies that also $\mathbb{P}(\|a_p\sqrt{n} \ \klaus{\HL}\|\leq C)>1-\varepsilon$, i.e $\klaus{\HL}=O_{\IP}(a_p^{-1}n^{-1/2})$. The approximation (A3) in~\cite{Mottonen2010} gives that for every $\delta\in (0,1)$
\[
\frac{a_p\sqrt{n}}{C}
\una{d_{2n}}\left(\frac{C}{a_p\sqrt{n}} \bo v\right)
\ge
-\sqrt{n}U_1+ C U_2 - C_1(C^{1+\delta}/n^{\delta/2}) U_3
\]
where
$
U_1={\binom{n}{2}}^{-1}\sum_{i<j} \textbf{u}_{i,j}^{\top}\textbf{v},\ \
U_2={\binom{n}{2}}^{-1}\sum_{i<j} \frac 1{2(a_pr_{i,j})}\{1-(\textbf{v}^{\top}\textbf{u}_{i,j})^2\}
$,
$
U_3={\binom{n}{2}}^{-1}\sum_{i<j} \frac 1{(a_pr_{i,j})^{1+\delta}},
$
\Una{and $C_1$ is a universal constant (not depending on $n,p$ or $\delta$).} The random variables $U_1$, $U_2$ and $U_3$ are all univariate U-statistics with
a well-developed theory. As $\mathbb{E}(r_{i,j}^{-k})<\infty$ for $k\leq 4$, these U-statistics have finite expectation and variances.  Further
\[
C U_2 - C_1(C^{1+\delta}/n^{\delta/2}) U_3=
 C \mathbb{E}[U_2] - C_1(C^{1+\delta}/n^{\delta/2}) \mathbb{E}[U_3] + U
\]
where $U$ is a U-statistic with the mean zero, finite variance,  and $\sqrt{n}U=O_{\IP}(1)$.
The lower bound can be written as
\[
\frac{a_p\sqrt{n}}{C} 
\una{d_{2n}}(Ca_p^{-1}n^{-1/2}\textbf{v})\geq -\sqrt{n} U_1 +
 C \mathbb{E}[U_2] - C_1(C^{1+\delta}/n^{\delta/2}) \mathbb{E}[U_3] + o_{\IP}(1).
\]
\Una{Assumption~\ref{assumption:HD1} (a)
ensures that $\mathbb{E}(1-(\textbf{v}^\top\textbf{u}_{1,2})^2)\geq C_2$, for some universal constant $C_2>0$, further giving $\mathbb{E}(U_2)\geq {C_2}/{2}$.} 
On the other hand, Assumption~\ref{assumption:HD1} (b) gives that $\mathbb{E}(U_3)=O(1)$. We thus obtain the lower bound 
\[
\frac{a_p\sqrt{n}}{C} 
\una{d_{2n}}(Ca_p^{-1}n^{-1/2}\textbf{v})\geq -\sqrt{n} U_1 +
 \frac{\Una{C_2}}{2}C + o_{\IP}(1).
\]
As $\mathrm{tr}(\mathbb{E}(\textbf{u}_{1,2}\textbf{u}_{1,2}^\top))=1$, the eigenvalues of the covariance of $\textbf{u}_{1,2}$ are uniformly bounded from above by $1$. Thus,  $\sqrt{n}U_1=\Una{O_p(1)}$ further implying that for every $\varepsilon>0$, we can find $C>0$, so that for $n>0$ large enough 
\[
\IP\{\Una{-\sqrt{n} U_1 +
 {{C_2}}C/2 + o_{\IP}(1)}
> 0\} > 1-\epsilon,
\]
thus completing the proof.
\end{myproof}

\begin{myproof}{Theorem}{\ref{theorem:HD_HL2}}
 The approximation (B2) 
 in \cite{Mottonen2010} gives that for $\boldsymbol{\mu}=\klaus{\HL}\in\mathbb{R}^p$ and $\delta\in(0,1)$
\begin{align}\label{eq:jana_B2}
    &\binom{n}{2}^{-1}\sum_{i<j}\frac{\textbf{z}_{i,j}-\boldsymbol{\mu}}{\|\textbf{z}_{i,j}-\boldsymbol{\mu}\|}=\binom{n}{2}^{-1}\sum_{i<j}\left(\textbf{u}_{i,j}\Una{-} r_{i,j}^{-1}(\textbf{I}_p-\textbf{u}_{i,j}\textbf{u}_{i,j}^\top)\boldsymbol{\mu}\right)+\textbf{R}(\boldsymbol{\mu}),
\end{align}
where the norm of the remainder $\|\textbf{R}(\boldsymbol{\mu})\|\leq C\|\boldsymbol{\mu}\|^{1+\delta}\binom{n}{2}^{-1}\sum_{i<j} 
r_{i,j}^{-(1+\delta)}$,
for some universal constant $C>0$ (does not depend on $n,\,p$). \Una{Observing that $\displaystyle\binom{n}{2}^{-1}\sum_{i<j}\frac{\textbf{z}_{i,j}-\klaus{\HL}}{\|\textbf{z}_{i,j}-\klaus{\HL}\|}=\klaus{\bo 0}$ and substituting $\klaus{\HL}=n^{-1/2}a_p^{-1}{\HL}^*$ in~\eqref{eq:jana_B2}, we obtain}
\[
\sqrt{n}\textbf{T}- \textbf{A} {\HL}^* = \textbf{R}({\HL}^*)
\]
where
\[
\textbf{T}={\binom{n}{2}}^{-1}\sum_{i<j} \textbf{u}_{i,j}
\ \ \ \mbox{and}\ \ \
\textbf{A}={\binom{n}{2}}^{-1}\sum_{i<j} \left[ \frac 1 {a_pr_{i,j}}(\textbf{I}-  \textbf{u}_{i,j} \textbf{u}_{i,j}^{\top}) \right]
\]
and
\[
\|\bo R({\HL}^*)\| \le  C_2 n^{-\delta/2} |{\HL}^*|^{1+\delta} \left[{\binom{n}{2}}^{-1}   \sum_{i<j}(a_p r_{i,j})^{-(1+\delta)}\right] 
\stackrel{\IP}\longrightarrow 0
\]
\bigskip
Next show that
$\bo D({\HL}^*)={\binom{n}{2}}^{-1}\sum_{i<j} \left[ \frac 1 {a_pr_{i,j}}\bo u_{i,j} \bo u_{i,j}^{\top} \right]{\HL}^*=o_{\IP}(1)$. Write first $\una{\boldsymbol{\mu}=c\textbf{v}}$ where $0<c<C$ and $\bo v'\bo v=1$. Then,
under sphericity,
$\bo D(\boldsymbol{\mu})$ is a U-statistics with expected value $(c/p)\bo v$ and
(by triangular inequality)
\[
\|\bo D(\una{\boldsymbol{\mu}})\|\le C {\binom{n}{2}}^{-1}\sum_{i<j} \left[ \frac 1 {a_pr_{i,j}} |\una{\textbf{u}_{i,j}^{\top}\textbf{v}}| \right]
\]
and the upper limit is a U-statistic with  the expected value at most $C/\sqrt{p}$
and the variance $O(1/(pn))$.  Therefore the upper limit
converges in probability to zero as $p\to\infty$. Finally, as ${\HL}^*=O_{\IP}(1)$,
also $\bo D({\HL}^*)\stackrel{\IP}\longrightarrow 0$.
This is seen as
\[
\IP(\|\bo D({\HL}^*)\|\le \epsilon)\le   \IP(\|{\HL}^*\|\le C)\cdot \IP(\|\bo D({\HL}^*)\|\le \epsilon \ |\ \|{\HL}^*\|\le C )+ \IP(\|{\HL}^*\|> C)
\]

This gives us the second approximation (under sphericity)
\[
\sqrt{n} \ \bo T-\left[{\binom{n}{2}}^{-1}\sum_{i<j}  \frac 1 {a_pr_{i,j}}\right]  \hat{\boldsymbol{\mu}}_{HL}^* = o_{\IP}(1)
\]
Finally as $\displaystyle{\binom{n}{2}}^{-1}\sum_{i<j} (a_pr_{i,j})^{-1}=1+(1/\sqrt{n})O_{\IP}(1)$, we have
the third approximation for ${\HL}^*=a_p \ \sqrt{n} \ \klaus{\HL}$:
\[
 {\HL}^*=\sqrt{n}  \una{\textbf{T}} +o_{\IP}(1). 
\]
\end{myproof}

\end{document}